\documentclass[11pt]{article}
\usepackage{color}
\usepackage{graphicx}
\usepackage{amsfonts}
\usepackage{theorem}
\newtheorem{Lem}{Lemma}[section]
\newtheorem{Def}[Lem]{Definition}
\newtheorem{The}[Lem]{Theorem}
\newtheorem{Prop}[Lem]{Proposition}
\newtheorem{Cor}[Lem]{Corollary}

\newtheorem{Cex}[Lem]{Counter-example}

\newtheorem{Rem}[Lem]{Remark}

\newcommand{\qed}{\hbox{\rule{6pt}{6pt}}}

\begin{document}
\title{On bounds of Tsallis relative entropy and an inequality for generalized skew information}
\author{Shigeru Furuichi$^1$\footnote{E-mail:furuichi@chs.nihon-u.ac.jp}\\
$^1${\small Department of Computer Science and System Analysis,}\\
{\small College of Humanities and Sciences, Nihon University,}\\
{\small 3-25-40, Sakurajyousui, Setagaya-ku, Tokyo, 156-8550, Japan}}
\date{}
\maketitle
{\bf Abstract.} Quantum entropy and skew information play important roles in quantum information science. 
They are defined by the trace of the positive operators so that the trace inequalities often have important roles to develop the mathematical theory in quantum information science.
In this paper, we study some properties for information quantities in quantum system through trace inequalities.
Especially, we give upper bounds and lower bounds of Tsallis relative entropy, which is a one-parameter extension of the relative entropy
in quantum system. In addition, we compare the known bounds and the new bounds, for both upper and lower bounds,  respectively.
We also give an inequality for generalized skew information by introducing a generalized correlation measure.
\vspace{3mm}

\vspace{3mm}


\vspace{3mm}

\section{Introduction}
Two important theorems in  quantum information theory have been proved in  \cite{Schu} and \cite{Hol,SW}.
In the paper \cite{Schu}, it was shown the relation
between quantum entropy (von Neumann entropy) \cite{Neu} and source coding theorem in quantum system (non-commutative system).
In the papers \cite{Hol,SW},  it was also shown the relation
between the Holevo bound, which is considered to be mutual information in quantum information theory, and coding theorem for the classical-quantum channel.
Especially, in these papers \cite{Hol,Schu,SW}, the role of von Neumann entropy and Holevo bound in quantum information were clarified 
so that quantum information theory has been progressed for around fifteen years \cite{NC,Petz}.  
Before such developments of quantum information, von Nuemann entropy and related entropies such as
 relative entropy \cite{Um} and mutual entropy \cite{Ohy} were studied in both direction from physics and mathematics \cite{OP,Weh}. 
The study for a generalization on information entropy in classical system (commutative system) has a long history and we have many literatures \cite{Csi,Ren,Tsa}. 
See also \cite{AD,ESS} and references therein.
As for one of the generalizations of entropy, we have studied the Tsallis entropy and the Tsallis relative entropy in quantum system 
\cite{FYK,FYK2,Furu1,Furu1_01,Furu_book,Furu1_02,Furu1_03}.
In the present paper, we study on the bounds of the Tsallis relative entropy which is a one-parameter extension of the Umegaki relative entropy.

As one of the mathematical studies on the topics about entropy theory, skew information \cite{WY1,WY2} and its concavity problem are famous. The concavity
problem for skew information was generalized by F.J.Dyson, and it was proven by E.H.Lieb in \cite{Lie}. 
It is also known that skew information presents the degree of noncommutativity between a certain quantum state represented by a density operator $\rho$
 (which is a positive operator with an unit trace) and an observable represented by self-adjoint operator $H$,  therefore an uncertainty relation using
skew information has been studied in \cite{FYK,Furuichi,Ko,Lin,Luo1,Luo2,Yanagi}.  
See \cite{Hei,Rob,Schr} for the original uncertainty relations which were represented by the trace inequalities and they show the uncertainty principle
which is a fundamental concept in quantum mechanical physics.
In the present paper, we give a trace inequality for a generalized skew information by introducing a generalized correlation measure.

\section{Upper bounds of Tsallis relative entropy}
Firstly we give some notation.
We denote the generalized exponential function by $\exp_{\nu}(x) \equiv (1 + \nu x)^{\frac{1}{\nu}}$ if $1+\nu x >0$,
otherwise it is undefined and its inverse function (generalized logarithmic function) by $\ln_{\nu}x \equiv \frac{x^{\nu} -1}{\nu}$,
for $\nu \in (0,1]$ and $x \geq 0$.
The functions $\exp_{\nu}(x) $ and $\ln_{\nu} x$ converge to $e^x$ and $\log x$ as $\nu \to 0$, respectively.
Note that the definition of the generalized logarithmic function $\ln_{\nu}(X)$ for a positive operator $X$ is well-defined.
In this paper, we define the generalized exponential function $\exp_{\nu}(X)$ for a positive operator $X$ by
$
\exp_{\nu}(X) \equiv \left(I + \nu X\right)^{\frac{1}{\nu}},
$
if $Tr\left[ \left(I + \nu X\right)^{\frac{1}{\nu}}\right] \in \mathbb{R}$.

The Tsallis relative entropy in quantum system (noncommutative system) 
is defined for the positive operators in the following manner.
For the study of entropies from mathematical viewpoint, the condition such as an unit trace is often relaxed.
See p.274 of \cite{Bha} or \cite{FYK,FYK2}. 
See also \cite{Tsa,Tsa1,Tsa2,Tsa3} for the original Tsallis entropy and its advances in statistical physics.

\begin{Def}   \label{Tsa_rela_ent_def}
The Tsallis relative entropy is defined by
$$
D_{\nu} (X\vert Y) \equiv \frac{\hbox{Tr}[X - X^{1-\nu}Y^{\nu}]}{\nu} = \hbox{Tr}[X^{1-\nu} (\ln_{\nu} X - \ln_{\nu} Y)]
$$
for positive operators $X, Y$ and $\nu \in (0,1]$.
\end{Def}

We have the following proposition, which gives an upper bound of the  Tsallis relative entropy. 
\begin{Prop} {\bf (\cite{FYK})}  \label{ub1}
For positive operators $X,Y$ and $\nu \in (0,1]$, the following inequality holds.
\begin{equation} \label{FYK2004_JMP}
D_{\nu}(X\vert Y) \leq -Tr\left[X \ln_{\nu}\left( X^{-1/2}YX^{-1/2}\right) \right].
\end{equation}
\end{Prop}

The further upper bound of the right hand side was given by T.Furuta in \cite{Furuta}, with the generalized Kantorovich constant:
$$
D_{\nu}(X\vert Y) \leq -Tr \left[X\ln_{\nu}\left( X^{-1/2}YX^{-1/2} \right)\right] \leq \left( \frac{1-K(\nu,h)}{\nu}\right) Tr[X]^{1-\nu} Tr[Y]^{\nu} + D_{\nu}(X\vert Y),
$$
where $X^{1/2}\ln_{\nu}\left(X^{-1/2}YX^{-1/2}\right)X^{1/2}$ is often called the Tsallis relative operator entropy for positive operators $X,Y$ 
and $K(\nu,h)$ is the generalized Kantorovich constant defined for $\nu \in \mathbb{R}$ and $h \in (0,\infty)$ with $h \neq 1$:
$$
K(\nu,h) \equiv \frac{(h^{\nu}-h)}{(\nu-1)(h-1)} \left( \frac{(\nu -1)}{\nu} \frac{h^{\nu}-1}{(h^{\nu}-h)} \right)^{\nu}.
$$ 
We note that $K(2,\frac{M}{m})=K(-1,\frac{M}{m})=\frac{(M+m)^2}{4mM}$, which is often called the Kantorovich constant. See \cite{Furuta_book,M-P} for details.
Since we have $K(0,h)=K(1,h)=1$ and $\frac{dK(\nu,h)}{d\nu}\vert_{\nu =0} = -\log S(h)$, 
the above inequalities recover the following inequalities (the first inequality below was originally given in \cite{HP}) :
$$
U(X\vert Y) \leq -Tr\left[X \log \left( X^{-1/2}YX^{-1/2} \right)\right] \leq Tr[X] \log S(h) + U(X\vert Y),
$$
in the limit $\nu \to 0$.
Where $U(X\vert Y)\equiv Tr[X\left(\log X -\log Y\right)]$ is the relative entropy introduced by Umegaki in \cite{Um},  $X^{1/2} \log\left( X^{-1/2}YX^{-1/2} \right) X^{1/2}$ is the
relative operator entropy introduced in \cite{BS,KF} and $S(h)\equiv \frac{h^{\frac{1}{h-1}}}{e \log h^{\frac{1}{h-1}}}$ is the Specht's ratio \cite{Sp}, 
where $h \equiv \frac{M_1M_2}{m_1m_2} > 1$ for $0< m_1 I \leq X \leq M_1 I$ and $0< m_2 I \leq Y \leq M_2 I$. (See \cite{Furuta} for details.)

We also have the  following proposition.
\begin{Prop} \label{ub2}
For positive operators $X,Y$ and $\nu \in (0,1]$, the following inequality holds.
\begin{equation}  \label{ineq_ub2}
D_{\nu}(X\vert Y) \leq \frac{Tr[(X-Y)_+]}{\nu},
\end{equation}
where $A_+ \equiv \frac{1}{2}\left( A+\vert A\vert \right)$ and $\vert A\vert \equiv (A^*A)^{1/2}$ for any operator $A$.
\end{Prop}

{\it Proof}:
It immediately follows from the trace inequality proven by K.M.R.Audenaert et.al. in \cite{Aud}:
\begin{equation}
Tr[A^sB^{1-s}] \geq \frac{1}{2} Tr[A+B-\vert A-B \vert]
\end{equation}
for positive operators $A,B$ and $s\in [0,1]$.

\hfill \qed

From the above two propositions, we have two different upper bounds for the Tsallis relative entropy.
It is quite natural to consider that we have the following inequality for positive operators $X$ and $Y$ and $\nu \in (0,1]$:
\begin{equation}  \label{comparison01_01}
-Tr\left[X \ln_{\nu}\left( X^{-1/2}YX^{-1/2}\right) \right] \leq \frac{Tr[(X-Y)_+]}{\nu},
\end{equation}
which is equivalent to
\begin{equation}  \label{comparison01_02}
Tr\left[ X\sharp_{\nu}Y \right] \geq \frac{1}{2} Tr[X+Y -\vert X-Y\vert ],
\end{equation}
where $X\sharp_{\nu}Y \equiv X^{1/2}(X^{-1/2}YX^{-1/2})^{\nu}X^{1/2}$ is  $\nu$-power mean. Here we note that we have
$$
Tr[X^{1-\nu}Y^{\nu}] \geq Tr\left[ X\sharp_{\nu}Y \right] 
$$
from the inequality (\ref{FYK2004_JMP}).  However we have the counter-example for the inequality (\ref{comparison01_02}) in the following.
We take $\nu=1/2$ and 
\[
X = \left( {\begin{array}{*{20}c}
   {10} & 7  \\
   7 & 5  \\
\end{array}} \right),Y = \left( {\begin{array}{*{20}c}
   {16} & 6  \\
   6 & 3  \\
\end{array}} \right).
\]
Then we have 
$$
Tr\left[ X\sharp_{\nu}Y \right] -\frac{1}{2} Tr[X+Y -\vert X-Y\vert ] \simeq -0.510619.
$$
Therefore the inequality (\ref{comparison01_01}) does not hold in general. 
That is, we can conclude that neither the inequality (\ref{FYK2004_JMP}) nor the inequality (\ref{ineq_ub2}) is uniformly better than the other.

\section{Lower bounds of Tsallis relative entropy}

We firstly note that the Tsallis relative entropy defined in Definition \ref{Tsa_rela_ent_def} is not always nonnegative. 
(If we impose on the condition such as an unit trace for two positive operators $X$ ad $Y$, then the Tsallis relative entropy has a nonnegativity.)
Therefore it is natural to have an interest in the lower bound of the Tsallis relative entropy defined for two positive operators. 
We have the following proposition, which gives a lower  bound of the Tsallis relative entropy. 

\begin{Prop} {\bf (\cite{FYK})}  \label{lb1}
For positive operators $X,Y$ and $\nu \in (0,1]$, the generalized Bogolivbov inequality holds.
\begin{equation}\label{lb1_ineq}
D_{\nu}(X\vert Y) \geq \frac{Tr[X]-(Tr[X])^{1-\nu}(Tr[Y])^{\nu}}{\nu}
\end{equation}
\end{Prop}

If we take the limit $\nu \to 0$, Proposition \ref{lb1} recovers the original Peierls-Bogoliubov inequality \cite{BLP,Hug}:
$$
U(X\vert Y) \geq Tr \left[ X \left(\log Tr[X] -\log Tr[Y]\right)  \right].
$$
We also easily find that the right hand side in the inequality (\ref{lb1_ineq}) is nonnegative, if we have the relation such that $Tr[X] \geq Tr[Y]$.
We also have the  following lower bound for the Tsallis relative entropy.
\begin{The}  \label{lb2}
For positive operators $X,Y$ and $\nu \in (0,1]$, if we have  $I \leq Y \leq X$, 
then we have the following inequality
\begin{equation}\label{lower_bound_rel_ent}
D_{\nu}(X\vert Y) \geq Tr\left[X^{1-\nu}\ln_{\nu} \left(Y^{-1/2}XY^{-1/2} \right)\right].
\end{equation}
\end{The}

It is notable that the condition $X \geq Y$ assures the nonnegativity of the right hand side in the inequality (\ref{lower_bound_rel_ent}).
To prove Theorem \ref{lb2}, we use the following lemmas.
\begin{Lem} {\bf (\cite{FM})}  \label{gen-GT}
For positive operators $X,Y$ and $\nu \in (0,1]$, we have
\begin{equation}\label{gen-GT-ineq}
Tr[  \exp_{\nu}(X+Y)  ] \leq Tr[\exp_{\nu}(X) \exp_{\nu}(Y)].
\end{equation}
\end{Lem}

Here we give a slightly different version of a variational expression for the Tsallis relative entropy.
 It can be proven by the similar way to Theorem 2.1 in \cite{Furu1}. 
However we give the proof for the convenience of readers as to be a self-contained article.

\begin{Lem}  \label{gen_rela_vari_exp}
For any $\nu \in (0,1]$ and any $d \in[0,\infty)$, we have the following relations.
\begin{itemize}
\item[(i)] If $A$ and $Y$ are positive operators, then we have
$$ d \ln_{\nu} \left(\frac{Tr[ \exp_{\nu}(A+\ln_{\nu}Y )]}{d} \right)= \max \left\{ Tr[X^{1-\nu}A] -D_{\nu} (X\vert Y) : X \geq 0, Tr[X]=d \right\}. $$
\item[(ii)] If $X$ and $B$ are positive operators with $Tr[X] = d$, then
$$  D_{\nu}(X\vert \exp_{\nu}(B))  =\max \left\{ Tr[X^{1-\nu} A] - d \ln_{\nu} \left(\frac{Tr[\exp_{\nu}(A+B)]}{d}\right) : A \geq 0 \right\}.$$
\end{itemize}
\end{Lem}

{\it Proof:}
For the case of $\nu=1$, it is trivial so that we assume $\nu \in (0,1)$.
Since we have $\lim_{x\to 0}x \ln_{\nu}\frac{a}{x}=0$ for $\nu \in (0,1)$,
it is also trivial for the case of $X=0$, thus we assume $X\neq 0$.
\begin{itemize}
\item[(1)]
We define
$$F_{\nu}(X) \equiv Tr[X^{1-\nu} A] -D_{\nu} (X\vert Y) $$
for a positive operator $X$ with $Tr[X] =d < \infty $. If we take the Schatten decomposition $X = \sum_{j=1}^{\infty} \mu_j E_j $, where all
$E_j$, $(j=1,2,\cdots ,\infty)$ are projections of rank one with $\sum_{j=1}^{\infty} E_j =I$ and $\mu_j \geq 0$, 
$(j=1,2,\cdots ,\infty)$ with $\sum_{j=1}^{\infty} \mu_j =d$,
then we rewrite
$$F_{\nu}\left(\sum_{j=1}^{\infty} \mu_j E_j\right) = \sum_{j=1}^{\infty}
\left\{ \mu_j^{1-\nu}Tr[E_j A] + \frac{1}{\nu} \mu_j^{1-\nu}  Tr[E_jY^{\nu}] - \frac{1}{\nu} \mu_j Tr[E_j]\right\}.$$
Then we have
$$
\frac{\partial^2}{\partial \mu_j^2} F_{\nu} \left( \sum_{j=1}^{\infty} \mu_j E_j\right)
=-\nu (1-\nu) \mu_j^{-\nu-1} Tr\left[E_j\left(A+\frac{1}{\nu}Y^{\nu}\right)\right] \leq 0,
$$
which means $F_{\nu}$ is concave function.
Thus we find $F_{\nu}(X)$ attains its maximum at a certain positive operator $X_0$ with $Tr[X_0] = d$.
Then for any self-adjoint operators $S$ with $Tr[S] =0$ (since for any $t\in \mathbb{R}$, $Tr[X_0 + t S] = d$ 
which is a condition on the positive operator
defined on the domain of the function $F_{\nu}$), there exists a positive operator $X_0$ such that
$$0 = \frac{d}{dt} F_{\nu}(X_0+tS) \vert _{t=0} = (1-\nu) Tr\left[S (X_0^{-\nu}A +\frac{1}{\nu} X_0^{-\nu}Y^{\nu} )\right],$$
 so that $X_0^{-\nu}A +\frac{1}{\nu} X_0^{-\nu}Y^{\nu} = c I$ for $c \in \mathbb{R}$. Thus we have
$$ X_0 = d\frac{\exp_{\nu}(A+\ln_{\nu}Y)}{Tr[\exp_{\nu}(A+\ln_{\nu}Y)]}.$$
by putting $c=\frac{1}{\nu}$ and satisfying the condition $Tr[X_0]=d$. By the formulae $\ln_{\nu}\frac{y}{x} = \ln_{\nu}y + y^{\nu} \ln_{\nu} \frac{1}{x}$
and $\ln_{\nu} \frac{1}{x} = -x^{-\nu} \ln_{\nu}x$, we have
\begin{eqnarray*}
F_{\nu}(X_0) &=& d^{1-\nu}\frac{Tr[\left\{\exp_{\nu}(A+\ln_{\nu}Y)\right\}^{1-\nu}A]}{Tr[\exp_{\nu}(A+\ln_{\nu}Y)]^{1-\nu}}\\
\hspace*{5mm}&& -d^{1-\nu} Tr\left[ \frac{\left\{\exp_{\nu}(A+\ln_{\nu}Y)\right\}^{1-\nu}}{Tr[\exp_{\nu}(A+\ln_{\nu}Y)]^{1-\nu}}
\left( \ln_{\nu}  \left(\frac{\exp_{\nu}(A+\ln_{\nu}Y)}{\frac{1}{d} Tr[\exp_{\nu}(A+\ln_{\nu}Y)]}\right) -\ln_{\nu}Y\right) \right]  \\
&=& d^{1-\nu}\frac{Tr[\left\{\exp_{\nu}(A+\ln_{\nu}Y)\right\}^{1-\nu}A]}{Tr[\exp_{\nu}(A+\ln_{\nu}Y)]^{1-\nu}}
-d^{1-\nu} Tr\left[ \frac{\left\{\exp_{\nu}(A+\ln_{\nu}Y)\right\}^{1-\nu}}{Tr[\exp_{\nu}(A+\ln_{\nu}Y)]^{1-\nu}}
\left\{ A+\ln_{\nu}Y \right.\right.\\
\hspace*{5mm}&&  \left.\left. +\left\{\exp_{\nu}(A+\ln_{\nu}Y)\right\}^{\nu}\ln_{\nu}\left(\frac{d}{Tr[\exp_{\nu}(A+\ln_{\nu}(A+\ln_{\nu}Y))]}\right) -\ln_{\nu}Y\right\} \right]  \\
&=&-d^{1-\nu} Tr[\exp_{\nu}(A+\ln_{\nu}Y)]^{\nu} \ln_{\nu} \left(\frac{d}{Tr[\exp_{\nu}(A+\ln_{\nu}Y)]}\right)\\
&=&d \ln_{\nu} \frac{Tr[\exp_{\nu}(A+\ln_{\nu}Y)]}{d}.
\end{eqnarray*}

\item[(2)] It follows from (1) that the functional
 $$g(A) \equiv d\ln_{\nu} \left(\frac{Tr[\exp_{\nu}(A+B)]}{d}\right),\quad d \equiv Tr[X] < \infty $$
defined on the set of all positive operator is convex, due to triangle inequality on $\max$.
Now let $A_0 = \ln_{\nu} X -B$, and define
$$G_{\nu}(A) \equiv Tr[X^{1-\nu}A] - d \ln_{\nu} \left(\frac{Tr[\exp_{\nu}(A+B)]}{d}\right), $$
which is concave on the set of all positive operator. Then for any self-adjoint operators $S$, there exists a positive operator $A_0$ such that
\begin{eqnarray*}
\frac{d}{dt}G_{\nu}(A_0+tS) \vert _{t=0} &=& Tr[X^{1-\nu}S] -d\left(\frac{Tr[X]}{d}\right)^{\nu-1} \frac{Tr[S(I+\nu \ln_{\nu}X)^{\frac{1-\nu}{\nu}}]}{d}\\
&=& Tr[X^{1-\nu}S]-Tr[SX^{1-\nu}]=0,
\end{eqnarray*}
using the formulae $\frac{d}{dx}\ln_{\nu}(x)=x^{\nu-1}$ and $\frac{d}{dx}\exp_{\nu}(x)=(1+\nu x)^{\frac{1}{\nu}-1}.$
Therefore $G_{\nu}(A)$ attaines the maximum
\begin{eqnarray*}
G_{\nu}(A_0) &=& Tr[X^{1-\nu}(\ln_{\nu}X -B)] -d\ln_{\nu}  \left(\frac{Tr[\exp_{\nu}(\ln_{\nu}X-B+B)]}{d}\right)\\
&=&Tr[X^{1-\nu}(\ln_{\nu}X -B)] -d \ln_{\nu} \left(\frac{Tr[X]}{d}\right)\\
&=& Tr[X^{1-\nu} (\ln_{\nu}X-\ln_{\nu}\exp_{\nu}(B))]\\
&=&D_{\nu}(X\vert \exp_{\nu}(B)).
\end{eqnarray*}
\end{itemize}

\hfill \qed

If we take $d=1$ and the limit $\nu \to 0$, then Lemma \ref{gen_rela_vari_exp} recovers Lemma 2.1 in \cite{HP} for positive operators $A$ and $B$.
It is notable that the original variational expressions for the Umegaki relative entropy proved by F.Hiai and D.Petz in  \cite{HP} holds for Hermitian matrices $A$ and $B$.
We are now in a position to prove Theorem \ref{lb2}.\\

{\it Proof of Theorem \ref{lb2}}:
Putting $B=\ln_{\nu}Y$ and $A=\ln_{\nu}Y^{-1/2}XY^{-1/2}$ in (ii) of Lemma \ref{gen_rela_vari_exp}
under the assumption of $I\leq Y\leq X$ which assures $A \geq 0$ and $B \geq 0$, and then using Lemma \ref{gen-GT},
we have
\begin{eqnarray*}
D_{\nu}(X\vert Y)  &=& D_{\nu}(X\vert \exp_{\nu}({\ln_{\nu}Y})) \\
&=& D_{\nu}(X\vert \exp_{\nu}(B)) \\
&\geq&  Tr[X^{1-\nu}A] - Tr[X] \ln_{\nu} \left(\frac{Tr[\exp_{\nu}(A+B)]}{Tr[X] }\right) \\
&\geq&  Tr[X^{1-\nu}A] - Tr[X] \ln_{\nu} \left(\frac{Tr[  \exp_{\nu}(A) \exp_{\nu}(B)    ]}{Tr[X] }\right) \\
&=& Tr[X^{1-\nu}\ln_{\nu} Y^{-1/2}XY^{-1/2} ]  -Tr[X] \ln_{\nu} \left(\frac{Tr[Y^{-1/2}XY^{-1/2} Y]}{Tr[X] }\right) \\
&=& Tr[X^{1-\nu}\ln_{\nu} Y^{-1/2}XY^{-1/2} ].  
\end{eqnarray*}

\hfill \qed

\begin{Rem}
\begin{itemize}
\item[(I)] The trace inequality (\ref{lower_bound_rel_ent}) is equivalent to the following trace inequality:
$$
Tr\left[ X^{1-\nu}\left\{ X^{\nu}-Y^{\nu}+I-\left(Y^{-1/2}XY^{-1/2}\right)^{\nu}     \right\}       \right] \geq  0.
$$
Therefore, if the following matrix inequality:
\begin{equation}  \label{op_lb}
X^{\nu}-Y^{\nu}+I-\left(Y^{-1/2}XY^{-1/2}\right)^{\nu}     \geq 0
\end{equation}
holds, then the trace inequality (\ref{lower_bound_rel_ent}) immediately holds. However the matrix inequality (\ref{op_lb}) does not hold in general, 
since we have the following counter-examples.
\begin{itemize}
\item[(i)] If we take $\nu=1$ and 
\[
X = \left( {\begin{array}{*{20}c}
   2 & 1  \\
   1 & 4  \\
\end{array}} \right),Y = \left( {\begin{array}{*{20}c}
   1 & 0  \\
   0 & 2  \\
\end{array}} \right),
\] 
satisfying the condition $I \leq Y \leq X$ (which is the assumption of Proposition \ref{lb2}), 
then one of the eigenvalues of the Hermitian matrix $X-Y+I-Y^{-1/2}XY^{-1/2}$ takes a negative value.
\item[(ii)] If  we take $\nu=1$ and  \[
X = \frac{1}{9}\left( {\begin{array}{*{20}c}
   2 & 1  \\
   1 & 5  \\
\end{array}} \right),Y = \frac{1}{3}\left( {\begin{array}{*{20}c}
   1 & 0  \\
   0 & 2  \\
\end{array}} \right),
\]
satisfying the condition $X\leq Y \leq I$, then one of the eigenvalues of the Hermitian matrix $X-Y+I-Y^{-1/2}XY^{-1/2}$ takes a negative value.
\end{itemize}

\item[(II)]Our next concern moves to the assumption of Proposition \ref{lb2}. 
We easily find that a counter-example for the trace inequality (\ref{lower_bound_rel_ent}), in the case that our  assumption $I\leq Y \leq X$ is not satisfied.
For example, if we take $\nu = 1$ and 
\[
X = \frac{1}{{15}}\left( {\begin{array}{*{20}c}
   {10} & { - 3}  \\
   { - 3} & {10}  \\
\end{array}} \right),Y = \frac{1}{{10}}\left( {\begin{array}{*{20}c}
   1 & 1  \\
   1 & 2  \\
\end{array}} \right),
\]
which does not satisfy the assumption $I\leq Y \leq X$ (but satisfy $0 < Y\leq X\leq I$),
then
$$
Tr[X-Y+I-Y^{-1/2}XY^{-1/2}] \simeq -20.9667.
$$
Thus the inequality (\ref{lower_bound_rel_ent}) does not hold in general for arbitrary positive operators $X$ and $Y$.
\end{itemize}
From (I) and (II), we may claim that Proposition \ref{lb2} is not a trivial result.
\end{Rem}

Closing this section, we give a comment on the comparison of two lower bounds for the Tsallis relative entropy.
Under the condition $I \leq Y \leq X$, we may have a conjecture such as
\begin{equation}  \label{comparison02_01}
Tr \left[X^{1-\nu} \ln_{\nu}\left(Y^{-1/2}XY^{-1/2}\right) \right] \geq \frac{Tr[X]-  (Tr[X])^{1-\nu}(Tr[Y])^{\nu}  }{\nu}
\end{equation}
for  positive operators $X$ and $Y$ and $\nu \in (0,1]$.
The inequality (\ref{comparison02_01}) is equivalent to the following inequality
\begin{equation}  \label{comparison02_02}
Tr\left[X^{1-\nu}\left(Y^{-1/2}XY^{-1/2}\right)^{\nu}\right] + (Tr[X])^{1-\nu}(Tr[Y])^{\nu}  \geq Tr[X^{1-\nu}] +Tr[X].
\end{equation}

Here we take two positive definite matrices \[
X = \left( {\begin{array}{*{20}c}
   {10} & 5  \\
   5 & 5  \\
\end{array}} \right),Y = \left( {\begin{array}{*{20}c}
   1 & 0  \\
   0 & 2  \\
\end{array}} \right),
\]
satisfying the condition $I \leq Y \leq X$.
Then for $\nu=0.1$, we have 
$$
Tr\left[X^{1-\nu}\left(Y^{-1/2}XY^{-1/2}\right)^{\nu}\right] + (Tr[X])^{1-\nu}(Tr[Y])^{\nu} -\left(Tr[X^{1-\nu}] +Tr[X]\right) \simeq 0.508133.
$$
For $\nu=0.9$, we also have
$$
Tr\left[X^{1-\nu}\left(Y^{-1/2}XY^{-1/2}\right)^{\nu}\right] + (Tr[X])^{1-\nu}(Tr[Y])^{\nu} -\left(Tr[X^{1-\nu}] +Tr[X]\right) \simeq -1.1696.
$$
Therefore the inequality (\ref{comparison02_01}) does not hold in general. 
That is, we can conclude that neither the inequality (\ref{lb1_ineq}) nor the inequality (\ref{lower_bound_rel_ent}) is uniformly better than the other.
This result supports that our Theorem \ref{lb2} is meaningful, in the sense of the comparison with Proposition \ref{lb1}.

\section{An inequality for a generalized skew information}
The uncertainty principle is a fundamental concept in quantum mechanical physics. 
It is represented by the famous Heisenberg uncertainty relation such as a trace inequality \cite{Hei}:
\begin{equation}   \label{HUL}
V_{\rho}(A) V_{\rho}(B) \geq \frac{1}{4}\vert Tr[\rho[A,B]]\vert^2 
\end{equation}
for a quantum state $\rho$ and two observables $A$ and $B$. Where 
the variance for a quantum state $\rho$ and an observable $H$ is defined by 
$V_{\rho}(H)\equiv Tr[\rho \left(H-Tr[\rho H]I\right)^2] = Tr[\rho H^2] -Tr[\rho H]^2$.
The further strong result was given by Schr\"odinger \cite{Schr}:
\begin{equation} \label{S_UL}
V_{\rho}(A) V_{\rho}(B)-\vert Re\left\{ Cov_{\rho}(A,B) \right\} \vert^2 \geq \frac{1}{4}\vert Tr[\rho[A,B]]\vert^2,  
\end{equation}
where the covariance is defined by $Cov_{\rho}(A,B) \equiv Tr[\rho \left(A-Tr[\rho A]I\right)\left(B-Tr[\rho B]I\right)].$
Due to its importance in quantum physics, the uncertainty relation has been studied by many researchers. 
Especially, some important results have been studied in the relation to the skew information representing a quantum uncertainty
from the viewpoints of quantum information science.
Here we firstly review about it. 
As it has been shown in \cite{Kos,Luo1,YFK}, we do not have the uncertainty relation type inequality for the Wigner-Yanase skew information \cite{WY1}:
\begin{equation}
I_{\rho}(H) \equiv \frac{1}{2} Tr\left[(i[\rho^{1/2},H_0])^2\right] =Tr[\rho H^2] -Tr[\rho^{1/2}H\rho^{1/2}H], \label{wy}   
\end{equation}
where $H_0\equiv H-Tr[\rho H]I$ for a  density operator $\rho$ and an observable $H$.
That is, the following trace inequality did not hold in general \cite{Kos,Luo1,YFK}:
\begin{equation}  \label{ul_wy}
I_{\rho}(A)I_{\rho}(B) \geq \frac{1}{4}\vert Tr[\rho[A,B]]\vert^2
\end{equation}
for a density operator $\rho$ and observables $A$ and $B$. Where $[X,Y]\equiv XY-YX$ is a commutator.
The counter example was given as follows.
\begin{Cex} {\bf (\cite{YFK})}  \label{cex01}
We take 
\[
\rho  = \frac{1}{4}\left( \begin{array}{l}
 3\,\,\,\,0 \\ 
 0\,\,\,\,1 \\ 
 \end{array} \right),A = \left( \begin{array}{l}
 \,\,0\,\,\,\,\,\,\,i \\ 
  - i\,\,\,\,\,\,0 \\ 
 \end{array} \right),B = \left( \begin{array}{l}
 0\,\,\,\,\,1 \\ 
 1\,\,\,\,\,\,0 \\ 
 \end{array} \right),
\]
then we have 
$I_{\rho}(A)I_{\rho}(B) = \left(1-\frac{\sqrt{3}}{2}\right)^2$ and $\frac{1}{4}\vert Tr[\rho[A,B]]\vert^2=\frac{1}{4}$.
Therefore the inequality (\ref{ul_wy}) does not hold in general.
\end{Cex}

As a one-parameter generalization, the Wigner-Yanase-Dyson skew information was defined by
\begin{eqnarray}  
I_{\rho,\alpha}(H) &\equiv& \frac{1}{2} Tr\left[(i[\rho^{\alpha},H_0])(i[\rho^{1-\alpha},H_0])\right] \label{wyd}    \\
 &=&Tr[\rho H^2] -Tr[\rho^{\alpha}H\rho^{1-\alpha}H], \quad \alpha\in[0,1] \nonumber 
\end{eqnarray}
where $H_0\equiv H-Tr[\rho H]I$ for a  density operator $\rho$ and an observable $H$.

In \cite{Luo2}, S.Luo introduced a new quantity such as
\begin{equation} \label{u1}
U_{\rho}(H) \equiv \sqrt{V_{\rho}(H)^2- \left(V_{\rho}(H)-I_{\rho}(H) \right)^2}
\end{equation}
for a  density operator $\rho$ and an observable $H$.
Then he succeeded to establish the uncertainty relation type inequality as follows:
\begin{equation}\label{ul_u1}
U_{\rho}(A)U_{\rho}(B) \geq \frac{1}{4}\vert  Tr[\rho[A,B]]\vert^2
\end{equation}
for a density operator $\rho$ and observables $A$ and $B$. He also introduced the quantity associated to
Wigner-Yanase skew information,
$$
J_{\rho}(H) \equiv \frac{1}{2} Tr\left[\left(i\left\{\rho^{1/2},H_0\right\}\right)^2     \right],
$$
where the anti-commutator is defined by $\{X,Y\} \equiv XY+YX$ for any operator $X$ and $Y$.
Then we have the relation $ U_{\rho}(H) = \sqrt{I_{\rho}(H) J_{\rho}(H)} $ and $V_{\rho}(H) = \frac{1}{2} \left( I_{\rho}(H) +J_{\rho}(H)\right)$.
Therefore, Luo's inequality (\ref{ul_u1}) refines Heisenberg's one (\ref{HUL}), since $U_{\rho}(H) \leq V_{\rho}(H)$.

K.Yanagi recently gave the generalization of the inequality (\ref{ul_u1}) as follows \cite{Yanagi}:
\begin{equation}\label{ul_u2}
U_{\rho,\alpha}(A)U_{\rho,\alpha}(B) \geq \alpha(1-\alpha)\vert  Tr[\rho[A,B]]\vert^2
\end{equation}
for $\alpha\in[0,1]$, a density operator $\rho$ and observables $A$ and $B$.
Where $U_{\rho,\alpha}(H)$ was defined by
\begin{equation} \label{u2}
U_{\rho,\alpha}(H) \equiv \sqrt{V_{\rho}(H)^2 - \left(V_{\rho}(H)-I_{\rho,\alpha}(H) \right)^2}
\end{equation}
for a  density operator $\rho$ and an observable $H$. 

In addition, quite recently, we gave the Schr\"odinger uncertainty relation for mixed states in \cite{Furuichi}:
\begin{equation} \label{furuichi_ineq}
U_\rho(A)U_\rho(B)-|Re\left\{ Corr_{\rho}(A,B)\right\}|^2 \geq  \frac{1}{4}|Tr[\rho[A,B]]|^2,
\end{equation}
where the correlation measure is defined for arbitrary operators $X$ and $Y$ by
$$
Corr_{\rho}(X,Y) \equiv Tr[\rho X^*Y] -Tr[\rho^{1/2} X^*\rho^{1/2}Y].
$$

On the other hand, S.Luo showed the trace inequality representing the relation between the original Wigner-Yanase skew information and the correlation measure
in \cite{Luo2}:
\begin{equation} \label{ineq_ul_i}
I_{\rho}(A) I_{\rho}(B) \geq \left|  Re \left\{Corr_{\rho}(A,B) \right\} \right|^2
\end{equation}
for a density operator $\rho$ and two observables $A$ and $B$.


It is remarkable that, if a quantum state (density operator) $\rho$ is a pure state (i.e., $\rho^2=\rho$), then the inequality (\ref{ineq_ul_i}) recovers
$$V_{\rho}(A)V_{\rho}(B) \geq  \frac{1}{4} \vert Tr[\rho \left\{A,B\right\}]\vert^2,$$
since $I_{\rho}(H)=V_{\rho}(H)$ and $Corr_{\rho}(A,B)=Cov_{\rho}(A,B)$ if $\rho$ is a pure state.
Where the covariance is defined by $Cov_{\rho}(A,B)\equiv Tr[\rho AB]-Tr[\rho A]Tr[\rho B]$.

In addition, defining a one-parameter extended correlation measure for $\alpha \in [0,1]$ and arbitrary operators $X$ and $Y$ by
$$
Corr_{\rho,\alpha}(X,Y) \equiv Tr[\rho X^*Y] -Tr[\rho^{\alpha}X^*\rho^{1-\alpha}Y].
$$
we have the following inequality:
\begin{equation} \label{ineq_ul_i_alpha}
I_{\rho,\alpha}(A) I_{\rho,\alpha}(B) \geq \left|  Re \left\{Corr_{\rho,\alpha}(A,B) \right\} \right|^2.
\end{equation}
for a density operator $\rho$, two observables $A$, $B$ and $\alpha\in[0,1]$, putting $\varepsilon = 0$ in Theorem III.4 of \cite{YFK}.

If we take $\alpha = \frac{1}{2}$, then the inequality (\ref{ineq_ul_i_alpha}) recovers the inequality (\ref{ineq_ul_i}).

To give the further generalized trace inequality, we give the following definition.
\begin{Def}
Let $f$ and $g$ be the operator monotone functions. 
Let $(f,g)$ be a monotonic pair. Where $(f,g)$ is called a monotonic pair if $(f(a) - f(b)) (g(a) - g(b)) \geq 0$ for any $a, b \in D\subset \mathbb{R}$,
for two functions $f$ and $g$ on the domain $D \subset \mathbb{R}$.
For a density operator $\rho$ and an observable $H$,
we define $(f,g)$-skew information $I_{\rho,(f,g)}(H)$ by
\begin{equation}
I_{\rho,(f,g)} (H) \equiv \frac{1}{2} Tr \left[\left(i\left[f(\rho),H_0\right]\right)(i\left[g(\rho),H_0\right]) \right],
\end{equation}
where $H_0\equiv H-Tr[\rho H]I$.
For a density operator $\rho$ and any operators $X,Y$,
we also define $(f,g)$-correlation measure $Corr_{\rho,(f,g)}(X,Y)$ by
\begin{equation}   \label{f-corr-def}
Corr_{\rho,(f,g)}(X,Y) \equiv Tr\left[f(\rho)g(\rho)X^*Y\right]-Tr\left[f(\rho)X^*g(\rho)Y\right].
\end{equation}
\end{Def}

Then we can prove the following theorem.
\begin{The}  \label{ul_w}
For $(f,g)$-skew informations $I_{\rho,(f,g)}(A)$, $I_{\rho,(f,g)}(B)$ and $(f,g)$-correlation measure $Corr_{\rho,(f,g)}(A,B)$, we have
\begin{equation}  \label{ineq_ul_w}
I_{\rho,(f,g)}(A) I_{\rho,(f,g)}(B) \geq   \left| Re \left\{Corr_{\rho,(f,g)}(A,B) \right\} \right|^2.
\end{equation}
\end{The}

{\it Proof}:
For $Corr_{\rho,(f,g)}(X,Y)$, we have the following properties, that is, $Corr_{\rho,(f,g)}(X,Y)$ is a sesquilinear form and Hermitian and
$Corr_{\rho,(f,g)}(A,A)$ has the nonnegativity for a self-adjoint operator $A$, since $(f,g)$ is a monotonic pair \cite{Bourin, Fujii}.
Then for self-adjoint operators $A$ and $B$, we have for any $t\in \mathbb{R}$
\begin{eqnarray*}
0 &\leq& Corr_{\rho,(f,g)}(tA+B,tA+B) \\
&=& Tr[f(\rho)g(\rho) (tA+B)(tA+B)]-Tr[f(\rho)(tA+B)g(\rho)(tA+B)] \\
&=& \left( Tr[f(\rho)g(\rho) A^2] -Tr[f(\rho)Ag(\rho)A]\right)  t^2 +2 Re\left( Tr[f(\rho)g(\rho) AB]-Tr[f(\rho)Ag(\rho)B] \right)t \\
&& \hspace*{5mm} + \left( Tr[f(\rho)g(\rho) B^2] -Tr[f(\rho)Bg(\rho)B]\right).
\end{eqnarray*}
Thus we have 
\begin{eqnarray*}
&&\left| Re\left\{ Tr[f(\rho)g(\rho) AB]-Tr[f(\rho)Ag(\rho)B] \right\}\right|^2 \\
&& \leq \left( Tr[f(\rho)g(\rho) A^2] -Tr[f(\rho)Ag(\rho)A]\right)  
 \left( Tr[f(\rho) g(\rho) B^2] -Tr[f(\rho) Bg(\rho) B]\right),
\end{eqnarray*}
which implies the theorem, since we have
$$
I_{\rho,(f,g)} (H) = Tr[f(\rho)g(\rho)H^2] -Tr[f(\rho)Hg(\rho)H].
$$

\hfill \qed

Theorem \ref{ul_w} recovers the inequality (\ref{ineq_ul_i_alpha}), putting $f(x)=x^{\alpha}$ and $g(x)=x^{1-\alpha}$ for $\alpha \in [0,1]$.
We also have the following corollary.
\begin{Cor}  \label{ul_k}
For $\alpha \in [0,1]$, a density operator $\rho$ and two observables $A$, $B$, we have
\begin{equation}  \label{ineq_ul_k}
K_{\rho,\alpha}(A) K_{\rho,\alpha}(B) \geq   \left| Re\left\{ Corr^{(K)}_{\rho,\alpha}(A,B) \right\} \right|^2,
\end{equation}
where a one-parameter extended correlation measure is  defined by
$$Corr^{(K)}_{\rho,\alpha}(X,Y) \equiv Tr\left[  \left(\frac{\rho^{\alpha}+\rho^{1-\alpha}}{2}\right)^2   X^*Y\right] 
-Tr\left[ \left(\frac{\rho^{\alpha}+\rho^{1-\alpha}}{2}\right) X^* \left(\frac{\rho^{\alpha}+\rho^{1-\alpha}}{2}  \right) Y\right]$$
for  $\alpha \in [0,1]$, a density operator $\rho$ and any operators $X$ and $Y$, 
and a one-parameter extended Wigner-Yanase skew information \cite{FYK3} is defined by $K_{\rho,\alpha}(H) \equiv Corr^{(K)}_{\rho,\alpha}(H,H) $
for  $\alpha \in [0,1]$, a density operator $\rho$ and a self-adjoint operator $H$.
\end{Cor}
{\it Proof}:
Put $f(x)=g(x)=\frac{x^{\alpha}+x^{1-\alpha}}{2}$ in Theorem \ref{ul_w}.

\hfill \qed

\section{Conclusion}
As we have seen, we have studied the properties of the fundamental information measure in quantum system, 
namely the Tsallis relative entropy and the generalized skew information through the trace inequality with the mathematical tools in matrix analysis.
Closing conclusion, we give the following quantity which generalizes the Tsallis relative entropy and the Wigner-Yanase-Dyson skew information.
For two positive operators $X$ and $Y$, and a self-adjoint operator $H$, we define
$$
L_t(X,Y;H) \equiv Tr[XH^2]-Tr[X^tHY^{1-t}H],\quad t \in [0,1].
$$
Then using $L_t(X,Y;H)$, the Wigner-Yanase-Dyson skew information can be rewritten by
$$I_{\alpha}(\rho,H)\equiv L_{\alpha}(\rho,\rho;H)=Tr[\rho H^2]-Tr[\rho^{\alpha} H\rho^{1-\alpha}H].$$
Also the Tsallis relative entropy can be rewritten by
$$D_{\nu}(X\vert Y)\equiv \frac{1}{\nu}L_{1-\nu}(X,Y;I) = \frac{Tr[X-X^{1-\nu}Y^{\nu}]}{\nu},\quad \nu \in (0,1].$$

It may be important to study the mathematical properties of the quantity $L_t(X,Y;H)$ in the future, 
since it covers both the Tsallis relative entropy (one-parameter extended relative entropy) and the Wigner-Yanase-Dyson skew information as special cases.

\section*{Acknowledgement}
The author was supported in part by the Japanese Ministry of Education, Science, Sports and Culture, Grant-in-Aid for
Encouragement of Young Scientists (B), 20740067

\end{document}